\documentclass{article}

\usepackage{arxiv}

\usepackage{orcidlink}

\usepackage[utf8]{inputenc} 
\usepackage[T1]{fontenc}    
\usepackage{hyperref}       
\usepackage{url}            
\usepackage{booktabs}       
\usepackage{amsfonts}       
\usepackage{nicefrac}       
\usepackage{microtype}      
\usepackage{lipsum}
\usepackage{graphicx}
\usepackage{authblk}
\usepackage{amsmath}
\usepackage{comment}
\usepackage{wrapfig}
\graphicspath{ {./images/} }

\title{Integer Programming Strategies for the Snake-in-the-Box Problem}

\author[1]{Tim Kunt \orcidlink{0009-0006-5732-3208}}
\author[2,3]{Thorsten Koch \orcidlink{0000-0002-1967-0077}}
\author[1]{Hyejun Park}

\affil[1]{Digital Data and Information for Society, Science, and Culture,  \protect\\ Zuse Institute Berlin, 
Takustr. 7, 14195 Berlin, Germany}

\affil[2]{Software and Algorithms for Discrete Optimization, \protect\\Technische Universität Berlin,
Straße des 17. Juni 135, 10623 Berlin, Germany}

\affil[3]{Applied Optimization, Zuse Institute Berlin,\\
Takustr. 7, 14195 Berlin, Germany}

\begin{document}
\maketitle

\begin{abstract}
We present a selection of Integer Programming (IP) strategies for the Snake-in-the-Box problem, which asks for the longest induced path in the $d$-dimensional hypercube.
Current records are based on sophisticated, exhaustive searches, with optimality proven only for instances up to dimension 8, motivating our study of IP formulations. We address symmetries and potential strengthenings. 
With the smallest open instance having as few as 512 binary decision variables, the proposed IPs offer a challenging and intriguing combinatorial problem. In addition, we extend the problem to higher-dimensional grid graphs and provide solutions for unexplored and unsolved instances.
\end{abstract}

\section{Introduction}
\label{Introduction}
\subsection{Problem Definition}
\label{Problem Definition}
Consider the hypercube $Q_d$ of dimension $d$ and the associated graph with the vertex set $V = \{0,1\}^d$ and the edge set 
$E = \bigl\{ \{u, v\} \;\big|\; u, v \in V,\; d_H(u,v) = 1 \bigr\}$, connecting each pair of vertices with Hamming distance ($d_H$) equal to $1$.\\
The Snake-in-the-Box Problem asks for the longest induced path on $Q_d = (V, E)$. An induced path is a path that is also an induced subgraph, meaning two non-consecutive vertices of the path cannot be adjacent. 

Analogously, the Coil-in-the-Box Problem asks for the longest induced cycle. The two problems are often addressed alongside each other, yet there are a few intricacies that distinguish the possible algorithms and techniques to solve them. 
For the generalization of the problem to grid graphs consider a $n^d$ lattice with $n\geq2$, and the associated graph $Q_{n^d}$ with the vertex set $V = \{0,\dots,n-1\}^d$ and the edge set $E = \bigl\{ \{u, v\} \;\big|\; u, v \in V,\; ||u-v||_1 = 1 \bigr\}$. 
\subsection{History and Literature Review}
\label{History and Literature Review}
\cite{kautz1958unit} is commonly cited as the first source to introduce the problem and provides its motivation as error detection codes 
and for circuit fault detection. 
To the best of our knowledge, the earliest proof of maximality of the snake of length $26$ on $Q_{2^6}$ dates back to 1965 \cite{davies2006longest}, followed up by the maximal snake for $Q_{2^7}$ at length $50$ in 1994 \cite{potter1994using}. Finally, a snake of length $98$ on $Q_{2^8}$ was first shown by \cite{carlson2010phenotype} in 2010 and proven maximal by \cite{ostergaard2015exhaustive} in 2015. At the time of this writing, the best lower bound for $Q_{2^9}$ is $191$, but it has yet to be proven maximal. 
This slow incremental progress is a testament to the difficulty of the problem and has been made possible through the advances in both computational power and algorithmic development. \\
\sloppy
\fussy
\begin{wrapfigure}[16]{r}{0.45\textwidth}
    \vspace{-\baselineskip}
    \centering
    \includegraphics[width=0.4\textwidth]{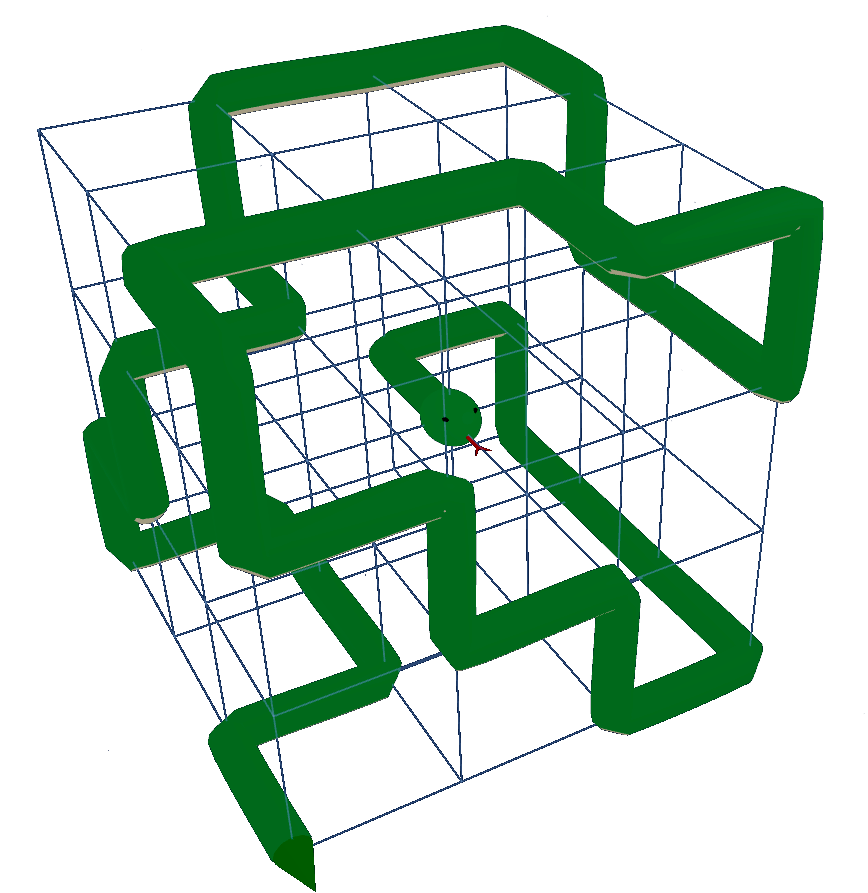}
    \caption{Longest snake on a $4\times 4\times4$ grid graph}
    \label{fig:1}
    \vspace{-10pt}
\end{wrapfigure}
\sloppy
Most recently, \cite{ace_2025_17344872}, \cite{ace_2025_18079790} and \cite{Echols} computed impressive improvements on lower bounds for snakes and coils. 
Just at the time of this writing in July 2026 \cite{orland2026census} further improved those bounds. An updated list of current records and certificates 
is available at \cite{AceWeb}.\\
\fussy
While this recent push has led to great improvements on the primal side of the problem (i.e., lower bounds), maximality of the proposed snakes and coils remains to be shown. 
A comparable IP approach was succesful in proving optimality for the generalized $n$-queens problem, and underlies many of the techniques used in the present work. 
For MIP formulations of the longest induced path and longest induced cycle problem on general graphs, we refer to 
\cite{melo2023mip} \cite{anaqreh2024exact} and \cite{bokler2020experimental}. A QUBO formulation for Snake-in-the-Box and Coil-in-the-Box has been proposed by \cite{fuidio2024qubo}.
\vspace{10pt}
\section{IP Formulations}
\label{IP Formulations}
Our aim is to prove maximality of snakes and coils on grid graphs of varying size and dimension. For that purpose we discuss and compare IP formulation and solution techniques tailored to the Snake-in-the-Box and Coil-in-the-Box problem. This problem class can be modelled as a cycle using subtour elimination routines, or with formulations that naturally avoid sub-cycles such as flow-based, walk and order-based formulations.

\subsection{Subtour Formulation}
\label{Subtour Formulation}
Following \cite{melo2023mip}, a promising and straightforward approach is to model snake and coil as a cycle with a subtour elimination routine. For this approach we may limit ourselves to the coil, as for the snake we may simply augment $Q_{n^d}$ by adding an auxiliary vertex to $Q_{n^d}^* := ( V^* := V \cup {a}, E^* := E \cup {\{av\}}_{v \in V})$, set $a$ as an active vertex, and find the longest cycle in $Q_{n^d}^*$. The solution corresponds to a maximal snake on $Q_{n^d}$, once we remove the two active edges containing $a$. 
Consider binary decision variables $y_e, x_v$ for edges and vertices. 
Let $\mathcal{C}$ denote all cycles in $Q_{n^d}$. 
\begin{align}
\max \quad & \sum_{e \in E} y_e \label{eq:tour_obj} \\[4pt]
\text{s.t.}\quad
& 1 + y_{(v_1,v_2)} \geq x_{v_1} + x_{v_2} 
  && \forall (v_1,v_2) \in E : v_1 \neq a,\; v_2 \neq a  \label{eq:tour_c1}, \\
& x_v \ge y_e
  && \forall v \in V,\; \forall e \in \delta(v) \label{eq:tour_c2}, \\
& \sum_{e \in \delta(v)} y_e = 2\,x_v
  && \forall v \in V \label{eq:tour_c3}, \\
  & \sum_{v\in c} x_v \leq |c|-1 && \forall c \in \mathcal{C}  \label{eq:cycles},\\
  & y_e, x_v \in \{0,1\}. \label{eq:int}
\end{align}
The objective function \eqref{eq:tour_obj} maximizes the number of active edges, equivalently vertices can be maximized over. Constraints \eqref{eq:tour_c1} enforce an induced cycle, in case of the augmented graph, an exception is made for edges containing the auxiliary node. Constraints \eqref{eq:tour_c2} set a vertex active if an adjacent edge is active. Constraints \eqref{eq:tour_c3} enforce that for each vertex exactly two adjacent edges are active if the vertex is active and none if not. Cycle constraints \eqref{eq:cycles} forbid solutions of multiple disjoint cycles. This can be formulated by either sub-tour elimination constraints or cut constraints, added dynamically. For a detailed discussion, we refer to \cite{melo2023mip} and \cite{bokler2020experimental}. \\
The above formulation can be achieved with only decision variables for edges, as in the model proposed by \cite{bokler2020experimental}. It can also be reformulated with decision variables for only the vertices.\\
The root relaxation for the base IP is weak, prompting our search of additional valid inequalities and exploitation of symmetries in the following section.\\
\subsection{Alternative Formulations}
Regarding alternative IP formulations, a walk-, flow-, and a potential order-based formulation all eliminate the need for a subtour elimination routine, but produce substantially larger models. During our initial experiments, the proposed subtour-elimination model outperformed all other variants, even with additional strengthenings, which aligns with the observations of \cite{bokler2020experimental}. 

\subsection{Additional valid inequalities and IP techniques}

\subsubsection{Bounds and Strengthenings}

A number of upper bounds are known and have been improved on, 
we refer to \cite{emelyanov2000maximal} and the references therein. 
Similarly good candidates for lower bounds are known and can be used to warm-start the IP. The difficulty of the IP model lies in closing the remaining dual gap.\\
\\
For every 4-cycle $C$ in $Q_{n^d}$ (for $n=2$ these are exactly the faces of the hypercube), we can further improve the LP relaxation. Unless $n=d=2$ and for $d\geq2$ it holds that 
\begin{align}
    \sum_{v \in C} x_v \leq 3 \quad \forall C \subset Q_{n^d},\\
    \sum_{e \in C} y_e \leq 2 \quad \forall C \subset Q_{n^d}.
\end{align}

In the trivial case $n=d=2$, this would cut off the optimal solution for the coil, which is the cycle around all four edges.\\
\\
Every grid graph is bipartite. For $v = (v_1, \dots, v_n) \in V$ consider
\begin{align}
    \pi(v) = \sum_{i=1}^{n} v_i \pmod{2}.
\end{align}
We call a vertex \textit{odd} if $\pi(v) = 1$ and \textit{even} if $\pi(v) = 0$ and write
\begin{align}
      V_{\mathrm{odd}}  &= \{ v \in V : \pi(v) = 1 \}, \\
      V_{\mathrm{even}} &= \{ v \in V : \pi(v) = 0 \}.
\end{align}
As the coil alternates between odd and even vertices, it holds that
\begin{align}
    \sum_{v \in V_{\mathrm{odd}}} x_v= \sum_{w \in V_{\mathrm{even}}} x_w. 
\end{align}
This constraint can be adapted to the snake, resulting in a weaker constraint. Following the same reasoning, for walk and order-based formulations, odd and even vertices can be constrained to odd and even steps.\\
\\
Consider $D$ the set of dimensions of $Q_{n^d}$ and $E_d$ the set of edges in each dimension direction, such that $(\{E_d\}_{d\in D}$ partitions $E$. A coil has an even number of edges in each dimension direction. In particular we may write
\begin{align}
& \sum_{e \in E_d} y_e = 2\,k_d
  && k \in \mathbb{Z}_{\geq0}, \forall d \in D.
\end{align}
The above constraints can be implemented by setting the $k_d$ as auxiliary integer variables of the model, which force the LHS to be even.\\
\\
For the coil, the IP benefits further from the fact that it has even length. Given a probable maximal candidate of length $l$, one can exchange the objective \eqref{eq:tour_obj} with the constraint $\sum_{e \in E} y_e \geq l+2$ (instead of just $l+1$). Infeasibility of the resulting IP proves that $l$ is the maximal length.

\subsubsection{Symmetries}

On $Q_{2^d}$, $d\geq3$ we may fix the first $4$ vertices of a snake or coil w.l.o.g. as 
\vspace{-2pt}
\begin{align*}
    (0,\dots,0,0,0,0), 
    (0,\dots,0,0,0,1), 
    (0,\dots,0,0,1,1), 
    (0,\dots,0,1,1,1).  
\end{align*} 
We can generalize this argument (c.f. canonical augmentation \cite{mckay1998isomorph}) to forbid multiple isomorphic snakes. We only allow the snake to enter a previously unvisited dimension if all other dimensions in lexicographic order have already been visited. 
Consider a snake as a walk on $Q_{2^d}$ starting at time $t=0$, then visiting a vertex at each time $t$ until some $\tau\leq T$. Define the binary decision variables $x^t_v$ active if and only if vertex $v$ is visited at time $t$. 
Let $V_d \subset V$ be the set of vertices for which the $d$-th entries are $1$. 
We get the following additional set of constraints for 
walk- and order-based formulations:
\begin{align}
    \sum_{i \in V_{d}}^{} 
    x^\tau_i \leq \sum_{j \in V_{d-1}}^{} \sum_{t=0}^{\tau-1} x^t_j & \quad \quad \forall \tau \in \{ 1, \dots , T\}, \forall d \label{eq:canonical}
\end{align}
The idea here is to enforce an order in which the snake or coil first enters each of the $d$ dimensions, which collapses isomorphic solutions and thereby massively shrinks the solution space of the IP. This does not cut off any valid solutions, as for any valid solution there exists an isomorphic solution that satisfies the above constraints. \\
There is a trade-off between breaking the symmetry of $Q_{n^d}$ in the formulation and exploiting a formulation- and solution-symmetric problem in the solving process. 
As a consequence of symmetry-breaking constraints, the IP is no longer formulation symmetric, which may limit the solver's own ability to exploit those symmetries. For that reason we additionally implemented orbital branching \cite{ostrowski2008orbital} in SCIP 10.0 \cite{hojny2025scip}, exploiting the automorphism group of the hypercube with a custom branch rule. We ran experiments on a set of instances using models implemented in both SCIP 10.0 and Gurobi 13.0.2 and with different formulations (subtour, walk-, and flow-based) and the discussed additional strengthenings. So far our computational experiments have not been entirely conclusive on which specific formulation is most suitable for breaking larger instances, considering all possible modelling choices. A custom branching rule in SCIP may break more symmetries than the solver finds by itself but could be outweighed by Gurobi's general advantage over SCIP. 
We observed that a dynamic approach using a subtour-elimination routine outperforms walk- and flow-based models, which aligns with the findings of \cite{bokler2020experimental}.

\section{Computational Results}

In the following tables, we compile the currently proven maximal snakes and coils on $n^d$-grid graphs \cite{OEIS}\cite{beluhov2023snake}. The special case $n=2$ is what is commonly referred to as Snake-in-the-Box. We visualized instances for $d=2$ and $d=3$ on a dedicated  \href{https://snake-in-the-box.vercel.app/}{website}\footnote{https://snake-in-the-box.vercel.app/}.\\
Results for $d\geq3, n\geq3$ are previously unknown maximum numbers of edges for snakes and coils that we have computed\footnote{experiments ran on Intel Xeon CPU E3-1270 v6 3.80GHz, Gurobi version 13.0.2, time limit of 1h} using the subtour-elimination routine with selected strengthenings. We were further able to reproduce the proofs for $n=2$ up to $d=7$.
\begin{table}[!htb]
    \centering

    \caption{maximal snakes on $Q_{n^d}$}
    \small
    \begin{tabular}{r|rrrrrrrrrr}
        $d/n$ &2&3&4&5&6&7&8&9&10&11  \\ \hline
        1 &   1& 2& 3& 4& 5& 6& 7& 8& 9&10  \\
        2 &   2& 6&10&16&23&32&41&52&63&76  \\
        3 &   4&15&35&66&  &  &  &  &  &    \\
        4 &   7&40&  &  &  &  &  &  &  &  \\
        5 &  13&  &  &  &  &  &  &  &  &  \\
        6 &  26&  &  &  &  &  &  &  &  &  \\
        7 &  50&  &  &  &  &  &  &  &  &  \\  
        8 &  98&  &  &  &  &  &  &  &  & \\  
    \end{tabular}

    \vspace{0.5cm}

    \caption{maximal coils on $Q_{n^d}$}
    \small
    \begin{tabular}{r|rrrrrrrrrr}
        $d/n$ &2&3&4&5&6&7&8&9&10&11  \\ \hline
        1 &   0& 0& 0& 0& 0& 0& 0& 0& 0& 0  \\
        2 &   4& 8&12&16&20&32&40&50&62&76   \\
        3 &   6&16&34&66&  &  &  &  &  &    \\
        4 &   8&40&  &  &  &  &  &  &  &  \\
        5 &  14&  &  &  &  &  &  &  &  &  \\
        6 &  26&  &  &  &  &  &  &  &  &  \\
        7 &  48&  &  &  &  &  &  &  &  &  \\  
        8 &  96&  &  &  &  &  &  &  &  & \\  
    \end{tabular}
\end{table}

\subsubsection*{Acknowledgements}
This work is co-funded by the European Union (European Regional Development Fund EFRE, Fund No. STIIV-001).

\bibliographystyle{alpha}
\bibliography{main}

@article{kautz1958unit,
  title={Unit-distance error-checking codes},
  author={Kautz, William H},
  journal={IRE Transactions on Electronic Computers},
  number={2},
  pages={179--180},
  year={1958},
  publisher={IEEE}
}

@article{davies2006longest,
  title={Longest Separated Paths and Loops in an N Cube},
  author={Davies, Donald W},
  journal={IEEE Transactions on Electronic Computers},
  number={2},
  pages={261--261},
  year={1965},
  publisher={IEEE}
}

@inproceedings{potter1994using,
  title={Using the Genetic Algorithm to Find Snake-in-the-Box Codes.},
  author={Potter, Walter D and Robinson, Robert W and Miller, John A and Kochut, Krys J and Redys, DZ},
  booktitle={IEA/AIE},
  pages={421--426},
  year={1994}
}

@inproceedings{carlson2010phenotype,
  title={Phenotype feedback genetic algorithm operators for heuristic encoding of snakes within hypercubes},
  author={Carlson, Benjamin P and Hougen, Dean F},
  booktitle={Proceedings of the 12th annual conference on Genetic and evolutionary computation},
  pages={791--798},
  year={2010}
}

@article{ostergaard2015exhaustive,
  title={Exhaustive search for snake-in-the-box codes},
  author={{\"O}sterg{\aa}rd, Patric RJ and Pettersson, Ville H},
  journal={Graphs and Combinatorics},
  volume={31},
  number={4},
  pages={1019--1028},
  year={2015},
  publisher={Springer}
}

@misc{ace_2025_17344872,
  author       = {Ace, Thomas},
  title        = {A New Lower Bound for Snake-in-the-Box in a
                   10-Dimensional Hypercube
                  },
  month        = oct,
  year         = 2025,
  publisher    = {Zenodo},
  doi          = {10.5281/zenodo.17344872},
  url          = {https://doi.org/10.5281/zenodo.17344872},
}

@misc{ace_2025_18079790,
  author       = {Ace, Thomas},
  title        = {New Lower Bounds for Snake-in-the-Box in 10-, 11-,
                   12-, and 13-dimensional Hypercubes
                  },
  month        = dec,
  year         = 2025,
  publisher    = {Zenodo},
  doi          = {10.5281/zenodo.18079790},
  url          = {https://doi.org/10.5281/zenodo.18079790},
}

@misc{Echols,
  author  = {Echols, William},
  title   = {Improving Snakes and Coils with Further Windowed Rerouting},
  year    = {2026},
  howpublished = {\url{https://williamechols.com}},
  note = {Last accessed July 2026}
}

@misc{AceWeb,
  title        = {Snake-in-the-Box lower bounds},
  author       = {Ace, Thomas},
  year         = 2026,
  howpublished = {\url{https://www.minortriad.com/snake/}},
  note         = {Last accessed July 2026}
}

@article{orland2026census,
  title={A Census of New Snake-in-the-Box Records},
  author={Orland, Paul and Fagan, Lucas and Tarquini, Michele and Passaro, Davide and Manko, Maksymilian and Heyes, Elli and Gruen, Angus and Butbaia, Giorgi and Tan, Justin and Gukov, Sergei},
  journal={arXiv preprint arXiv:2607.15270},
  year={2026}
}

@article{melo2023mip,
  title={{MIP} formulations for induced graph optimization problems: a tutorial},
  author={Melo, Rafael A and Ribeiro, Celso C},
  journal={International Transactions in Operational Research},
  volume={30},
  number={6},
  pages={3159--3200},
  year={2023},
  publisher={Wiley Online Library}
}

@article{anaqreh2024exact,
  title={Exact methods for the longest induced cycle problem},
  author={Anaqreh, Ahmad Turki and G-T{\'o}th, Bogl{\'a}rka and Vink{\'o}, Tam{\'a}s},
  journal={Croatian Operational Research Review},
  volume={15},
  number={2},
  pages={199--212},
  year={2024},
  publisher={Hrvatsko dru{\v{s}}tvo za operacijska istra{\v{z}}ivanja}
}

@inproceedings{bokler2020experimental,
  title={An experimental study of {ILP} formulations for the longest induced path problem},
  author={B{\"o}kler, Fritz and Chimani, Markus and Wagner, Mirko H and Wiedera, Tilo},
  booktitle={International Symposium on Combinatorial Optimization},
  pages={89--101},
  year={2020},
  organization={Springer}
}

@article{mckay1998isomorph,
  title={Isomorph-free exhaustive generation},
  author={McKay, Brendan D},
  journal={Journal of Algorithms},
  volume={26},
  number={2},
  pages={306--324},
  year={1998},
  publisher={Elsevier}
}

@techreport{ostrowski2008orbital,
  title={Orbital branching},
  author={Ostrowski, James and Linderoth, Jeff and Rossi, Fabrizio and Smriglio, Stefano},
  year={2008},
  institution={University of Wisconsin-Madison Department of Computer Sciences}
}

@article{hojny2025scip,
  title={The SCIP Optimization Suite 10.0},
  author={Hojny, Christopher and Besan{\c{c}}on, Mathieu and Bestuzheva, Ksenia and Borst, Sander and Dion{\'\i}sio, Jo{\~a}o and Ehls, Johannes and Eifler, Leon and Ghannam, Mohammed and Gleixner, Ambros and G{\"o}{\ss}, Adrian and others},
  journal={arXiv preprint arXiv:2511.18580},
  year={2025}
}

@article{emelyanov2000maximal,
  title={On the maximal length of a snake in hypercubes of small dimension},
  author={Emelyanov, Pavel G and Lukito, Agung},
  journal={Discrete Mathematics},
  volume={218},
  number={1-3},
  pages={51--59},
  year={2000},
  publisher={Elsevier}
}

@article{beluhov2023snake,
  title={Snake Paths in King and Knight Graphs},
  author={Beluhov, Nikolai},
  journal={arXiv preprint arXiv:2301.01152},
  year={2023}
}

@misc{OEIS,
  title        = {A099155, A000937, A357357},
  author       = {The On-Line Encyclopedia of Integer Sequences},
  howpublished = {\url{https://oeis.org}},
  note         = {Last accessed July 2026}
}

@article{fuidio2024qubo,
  title={{QUBO} formulation for the Snake-in-the-box and Coil-in-the-box problems},
  author={Fuidio, Federico and Canale, Eduardo and Sotelo, Rafael},
  journal={arXiv preprint arXiv:2409.04476},
  year={2024}
}

\end{document}